\documentclass[12pt]{amsart}

\usepackage{amsmath,amsthm,amssymb}

\font\teneufm=eufm10 \font\seveneufm=eufm7 \font\fiveeufm=eufm5
\newfam\frakturfam
\textfont\frakturfam=\teneufm \scriptfont\frakturfam=\seveneufm
\scriptscriptfont\frakturfam=\fiveeufm

\newtheorem{theor}{Theorem}
\newtheorem{co}{Corollary}

\def\bee{\begin{eqnarray}}
\def\bes{\begin{eqnarray*}}
\def\eee{\end{eqnarray}}
\def\ees{\end{eqnarray*}}

\def\a{\alpha}

\def\Proof{{\sl Proof.}\ }

\title{Identities of dual Leibniz algebras}

\begin{document}
\date{}
\maketitle

\begin{center}

{\bf Altyngul Naurazbekova}\footnote{Supported
by a grant of Kazakhstan,
Department of Mathematics, Eurasian National University,
 Astana, 010008, Kazakhstan,
e-mail: {\em altyngul.82@mail.ru}}
and
{\bf Ualbai Umirbaev}\footnote{Supported by NSA grant 3-00273 and by a grant of Kazakhstan, Department of Mathematics, Eurasian National University,
 Astana, 010008, Kazakhstan,
Department of Mathematics, Wayne State University,
Detroit, MI 48202, USA,
e-mail: {\em umirbaev@math.wayne.edu}}

\end{center}

\begin{abstract}
We prove that in characteristic $0$ any proper subvariety
of the variety of dual Leibniz algebras is nilpotent.
Consequently, the variety of dual Leibniz algebras
is Shpekhtian and has base rank 1.
\end{abstract}

\indent\hspace{6.5mm}

\indent\hspace{6.5mm}

\bigskip
\bigskip

\noindent
\begin{center}
{INTRODUCTION}
\end{center}
\hspace*{\parindent}
\bigskip

Recall that an algebra $\mathfrak{g}$ over an arbitrary field $k$ equipped with a bilinear operation $[-,-]$ is called {\em Leibniz} if it satisfies the (right) {\em Leibniz identity}
\bes
[x,[y,z]]=[[x,y],z]-[[x,z],y].
\ees
An algebra $A$ over $k$ is called {\em
dual Leibniz} if it satisfies the identity
\bee\label{f1}
(xy)z=x(zy+yz).
\eee
The Leibniz algebras form a Koszul operad in the sense of V. Ginzburg and M. Kapranov \cite{GK}.
Under the Koszul duality the operad of Lie algebras is dual to the operad of associative and commutative algebras.
The notion dual Leibniz algebra defined by J.-L. Loday \cite{Loday95} is precisely the dual operad of Leibniz algebras in this sense. Moreover, any dual algebra $A$ with respect to symmetrization
\bee\label{f2}
a\ast b= ab+ba
\eee
is an associative and commutative algebra \cite{Loday95}. This defines a functor
$$(dual \ \ Leibniz)\longrightarrow (Com)$$ between the categories of algebras,
which is dual to the inclusion functor $$(Lie)\longrightarrow (Leibniz).$$
Dual Leibniz algebras are called also Zinbiel (Leibniz is written in reverse order) algebras.

In this paper we study identities of dual Leibniz algebras. Let $\mathfrak{M}$ be a variety of linear algebras over a field $k$ and $k_\mathfrak{M}<x_1,x_2,\ldots,x_n>$ be the free algebra of $\mathfrak{M}$ in the variables $x_1,x_2,\ldots,x_n$. The least natural number $n$ such that the variety $Var(k_\mathfrak{M}<x_1,x_2,\ldots,x_n>)$ of algebras generated by  $k_\mathfrak{M}<x_1,x_2,\ldots,x_n>$ is equal to $\mathfrak{M}$ is called the {\em base rank} $rb(\mathfrak{M})$ of the variety $\mathfrak{M}$. If such a number does not exist then we say that $rb(\mathfrak{M})=\infty$.

It is obvious that the base rank of the variety of associative and commutative algebras is equal to $1$.
In 1952 A.I. Malcev \cite{Malcev52} proved that any associative algebra of countable dimension can be embedded in an associative algebra with two generators. The same result for Lie algebras was proved by A.I. Shirshov \cite{Shirshov58} in 1958. Consequently, the base ranks of the variety of associative algebras   and Lie algebras are equal to $2$. In 1977 I.P. Shestakov \cite{Shestakov77} proved that the base ranks of the variety of alternative and Malcev algebras are infinite (see, also \cite{Filippov80}).

A variety of algebras $\mathfrak{M}$ is said to be {\em Spechtian} if each of its subvarieties is defined by a finite system of identities. It is equivalent that the variety satisfies the descending chain condition for subvarieties.
The famous result by A.R. Kemer \cite{Kemer87} says that the variety of associative algebras is Shpekhtian in characteristic $0$. In the case of nonassociative algebras there are many partial results (see, for example \cite{Drensky00,Iltyakov92,Umirbaev84,Umirbaev85}).

A.S. Dzhumadildaev and K.M. Tulenbaev \cite{Dzhuma05}
proved an analog of the Nagata-Higman theorem \cite{Higman} for dual Leibniz algebras. In particular, they proved that any dual Leibniz algebra with bounded nil-index is nilpotent and every finite-dimensional dual Leibniz algebra over a field of characteristic $0$ is nilpotent.

J.-L. Loday proved \cite{Loday95} that the set of all nonassociative words with left arranged parenthesis compose a basis of free dual Leibniz algebras. So, there is a one-to-one correspondence between this basis of free dual Leibniz algebras and the set of all nonempty associative words. In this sense the properties of the variety of dual Leibniz algebras must be close to properties of the variety of associative algebras. This opinion contradicts to the results of this paper.

In this paper we prove that in characteristic $0$ the free dual Leibniz algebras in a single variable do not satisfy any nontrivial identity. We also prove that any proper subvariety of dual Leibniz algebras is nilpotent. Consequently, the variety of dual Leibniz algebras is Shpekhtian and has base rank 1.
It looks like there is some parallel between free dual Leibniz algebras and polynomial algebras which approves Loday's functor
$(dual \ \ Leibniz)\longrightarrow (Com)$.

The first author wishes to thank the Department of Mathematics of
Wayne State University in Detroit for the support while she was
working on this project.

\section{Homomorphisms into one generated free algebras}

\hspace*{\parindent}

Throughout this paper denote by $k$ an arbitrary fixed field of characteristic $0$.
Denote by $DL<x_1,x_2,\ldots,x_n,\ldots>$ the free dual algebra over $k$
in the variables $x_1,x_2,\ldots,x_n,\ldots$.

A linear basis of free dual Leibniz algebras is
given in \cite{Loday95}. The set of all nonempty nonassociative words with left arranged parenthesis
\bes
x_{i_1}(x_{i_2}\ldots(x_{i_{n-1}}x_{i_n})\ldots),  \ \ \ \ n\geq 1,
\ees
compose a basis of $DL<x_1,x_2,\ldots,x_n,\ldots>$.

In characteristic $0$ any identity is equivalent to the set of multilinear homogeneous identities \cite{KBKA}.
Then any nontrivial dual Leibniz identity can be written as
\bee\label{f3}
\sum_{\sigma \in S_n} \a_{\sigma}(x_{\sigma 1}(x_{\sigma 2}\ldots(x_{\sigma {n-1}}x_{\sigma n})\ldots))=0,
\eee
where $S_n$ is the $n$th symmetric group, $\a_{\sigma}\in k$

Let $A$ be an arbitrary dual Leibniz algebra and $a\in A$. We define $a^i$ by induction on $i$ as follows:
$a^1=a, a^{i+1}=a a^i$ for all $i\geq 1$.
It is proved in  \cite{Dzhuma05} that
\bee\label{f4}
a^ia^j=
\left(\begin{array}{cc}
i+j-1\\
j\\
\end{array}\right) a^{i+j}
\eee
for all $i,j\geq 1$.

\begin{theor}\label{t1}
The free dual Leibniz algebras in a single variable over a field of characteristic zero do not satisfy any nontrivial dual Leibniz identity.
\end{theor}
\Proof Denote by $A=DL<x>$ the free dual Leibniz algebra in the variable $x$.
Put $X_i=i! x^i$. By (\ref{f4}),
\bee\label{f5}
X_i X_j=\frac{i}{i+j}X_{i+j}.
\eee
Consider a homomorphism
\bes
\psi: DL<x_1,x_2,\ldots,x_n> \longrightarrow D<x>
\ees
such that $\psi(x_i)=X_{\lambda_i}$ for all $1\leq i\leq n$. Put also
\bes
P_n(\lambda_1,\lambda_2,\ldots,\lambda_n)=
\frac{\lambda_1}{\lambda_1+\lambda_2+\ldots+\lambda_n}\frac{\lambda_2}{\lambda_2+\ldots+\lambda_n}\ldots \frac{\lambda_{n-2}}{\lambda_{n-2}+\lambda_{n-1}+\lambda_n}\frac{\lambda_{n-1}}{\lambda_{n-1}+\lambda_n}.
\ees
Note that
\bes
P_n(\lambda_1,\lambda_2,\ldots,\lambda_n)=
\frac{\lambda_1}{\lambda_1+\lambda_2+\ldots+\lambda_n}P_{n-1}(\lambda_2,\ldots,\lambda_n).
\ees
Using this property and (\ref{f5}), it can be easily proved that
\bes
\psi(x_1(x_2\ldots(x_{n-1}x_n)\ldots))=P_n(\lambda_1,\lambda_2,\ldots,\lambda_n) X_{\lambda_1+\lambda_2+\ldots+\lambda_n}.
\ees
Put also
\bes
Q_n(\lambda_1,\lambda_2,\ldots,\lambda_n)=
\frac{1}{(\lambda_2+\ldots+\lambda_n)(\lambda_3+\ldots+\lambda_n)\ldots (\lambda_{n-2}+\lambda_{n-1}+\lambda_n)(\lambda_{n-1}+\lambda_n)\lambda_n}.
\ees
Then,
\bes
P_n(\lambda_1,\lambda_2,\ldots,\lambda_n)=
\frac{\lambda_1\lambda_2\ldots \lambda_n}{\lambda_1+\lambda_2+\ldots+\lambda_n}Q(\lambda_1,\lambda_2,\ldots,\lambda_n).
\ees

Now assume, in contrary, that $A$ satisfies a nontrivial multilinear identity of the form (\ref{f3}) with the least possible natural $n$. We may assume $\a_1=1$. Applying homomorphism $\psi$, from  (\ref{f3}) we get
\bes
\sum_{\sigma \in S_n} \a_{\sigma}P_n(\lambda_{\sigma 1},\lambda_{\sigma 2},\ldots,\lambda_{\sigma n})X_{\lambda_1+\lambda_2+\ldots+\lambda_n}=0.
\ees
Consequently,
\bes
\sum_{\sigma \in S_n} \a_{\sigma}P_n(\lambda_{\sigma 1},\lambda_{\sigma 2},\ldots,\lambda_{\sigma n})=0
\ees
and
\bes
\sum_{\sigma \in S_n} \a_{\sigma}Q_n(\lambda_{\sigma 1},\lambda_{\sigma 2},\ldots,\lambda_{\sigma n})=0.
\ees
The last equation can be written as
\bee\label{f6}
\sum_{\sigma \in S_n, \sigma n=n} \a_{\sigma}Q_n(\lambda_{\sigma 1},\lambda_{\sigma 2},\ldots,\lambda_{\sigma (n-1)},\lambda_n)
+
\sum_{\sigma \in S_n, \sigma n\neq n} \a_{\sigma}Q_n(\lambda_{\sigma 1},\lambda_{\sigma 2},\ldots,\lambda_{\sigma n})=0.
\eee
Thus, every system of positive integers $\lambda_1,\lambda_2,\ldots,\lambda_n$ satisfies equation (\ref{f6}).
Let's fix arbitrary positive integers $\lambda_1,\lambda_2,\ldots,\lambda_{n-1}$ and consider this equation  with respect to the variable $\lambda_n$. Multiplying the right hand side of the equation by common denominator, we get a polynomial with respect to $\lambda_n$. This polynomial is identically zero since it has infinite number of positive integer roots. Consequently, equation (\ref{f6}) holds for every $\lambda_n\in k$ if the right hand side just defined.

Note that the denominator of $Q_n(\lambda_{\sigma 1},\lambda_{\sigma 2},\ldots,\lambda_{\sigma n})$ is divisible by $\lambda_n$ if and only if $\sigma n=n$. Note also that by multiplying $Q_n(\lambda_1,\lambda_2,\ldots,\lambda_n)$ by $\lambda$ and substituting $\lambda_n=0$ we obtain $Q_{n-1}(\lambda_1,\lambda_2,\ldots,\lambda_{n-1})$.
Then, by multiplying equation (\ref{f6}) by $\lambda_n$ and substituting $\lambda_n=0$, we get
\bes
\sum_{\sigma \in S_n, \sigma n=n} \a_{\sigma}Q_{n-1}(\lambda_{\sigma 1},\lambda_{\sigma 2},\ldots,\lambda_{\sigma (n-1)})\\
=\sum_{\delta \in S_{n-1}} \a_{\delta}Q_{n-1}(\lambda_{\delta 1},\lambda_{\delta 2},\ldots,\lambda_{\delta (n-1)})
=0
\ees
for all positive integers $\lambda_1,\lambda_2,\ldots,\lambda_{n-1}$. Consequently,
\bes
\sum_{\delta \in S_{n-1}} \a_{\delta}P_{n-1}(\lambda_{\delta 1},\lambda_{\delta 2},\ldots,\lambda_{\delta (n-1)})
=0
\ees
and
\bes
\psi(\sum_{\delta \in S_{n-1}} \a_{\delta}(x_{\delta 1}(x_{\delta 2}\ldots(x_{\delta (n-2)}x_{\delta (n-1)})\ldots)))\\
=\sum_{\delta \in S_{n-1}} \a_{\delta}P_{n-1}(\lambda_{\delta 1},\lambda_{\delta 2},\ldots,\lambda_{\delta (n-1)})
X_{\lambda_1+\lambda_2+\ldots +\lambda_{n-1}}=0.
\ees
this means that algebra $A=DL<x>$ satisfies the identity
\bes
\sum_{\delta \in S_{n-1}} \a_{\delta}(x_{\delta 1}(x_{\delta 2}\ldots(x_{\delta (n-2)}x_{\delta (n-1)})\ldots))=0.
\ees
This identity is nontrivial  since $\a_1=1$. It contradicts the minimality of $n$. $\Box$

\section{Corollaries and applications}

\hspace*{\parindent}

\begin{co}\label{c1}
The base rank of the variety of dual Leibniz algebras is equal to one.
\end{co}
\Proof
By Theorem \ref{t1}, the variety of algebras $Var(DL<x>)$ generated by the free dual Leibniz algebra $DL<x>$ in the variable $x$ is defined only by identity (\ref{f1}). This means that the base rank of the variety of dual Leibniz algebras equals one.
$\Box$

It is easy to see that free dual Leibniz algebras in more than one variable cannot be embedded into one generated free dual Leibniz algebras.
\begin{theor}\label{t2}
Any proper subvariety of dual Leibniz algebras in characteristic zero is nilpotent.
\end{theor}
\Proof Let $\mathfrak{M}$ be an arbitrary proper subvariety of the variety of dual Leibniz algebras and
$B=\mathfrak{M}<y>$ be the free algebra of this variety in the variable $y$. Consider the homomorphism
$$\varphi : DL<x>\longrightarrow B$$
such that $\varphi(x)=y$. If $\varphi$ is isomorphism then $DL<x>\in \mathfrak{M}$, which is impossible by Theorem \ref{t1}. Consequently, $Ker(\varphi)\neq 0$. Then there exists a natural $n$ such that $x^n\in Ker(\varphi)$ since $\varphi$ is a homogeneous homomorphism of homogeneous algebras. This means that $\mathfrak{M}$
satisfies the identity
\bes
y^n=0,
\ees
i.e., the variety of algebras $\mathfrak{M}$ has nil index $n$. It was proved in \cite{Dzhuma05} that
nil algebras of bounded nil-index are nilpotent. Consequently, the variety $\mathfrak{M}$ is nilpotent. $\Box$

Note that every nilpotent variety of algebras is Shpekhtian, i.e., its every subvariety has finite basis of identities.
\begin{co}\label{c2}
The variety of dual Leibniz algebra in characteristic zero is Shpekhtian.
\end{co}

It is well known that any finite dimensional algebra satisfies an analogue of the standard identity (see definition, for example in \cite{Kemer87}). Then next corollary immediately follows from Theorem \ref{t2}.
\begin{co}\label{c3}
Every finite dimensional dual Leibniz algebra in characteristic zero is nilpotent.
\end{co}

This result was proved also in \cite{Dzhuma05}.

\end{document}